
\magnification=1000
\tolerance 1000
\vsize=21cm
\hsize=14cm
\baselineskip 18pt

\def\ska{\vskip 1truecm}
\def\ske{\vskip .5cm}

\def\half{{1\over 2}}
\def\ctrline#1{\line{\hss#1\hss}}

\def\a{\alpha}
\def\b{\beta}

\def\s{\sigma}
\def\D{\Delta}

\def\g{\gamma}

\def\ra{\rightarrow}
\def\st{\subset}
\def\half{{1\over 2}}
\def\nty{\infty}

\def\pr{^{\prime}}
\def\black{\penalty 10000\quad\penalty 10000}

\def\tf{\vbox{\vskip 2.5truecm}}

\newdimen\offdimen
\def\offset#1#2{\offdimen #1
   \noindent \hangindent \offdimen
   \hbox to \offdimen{#2\hfil}\ignorespaces}

\def\bk{\hfil\break}

\mathchardef\Dscr"244

\mathchardef\Kscr"24B  

\mathchardef\Rscr"252

\mathchardef\lscr"160

\def\corr{correspond}

\def\Hf{Hopf algebra}

\def\li{Lie algebra}
\def\mt{multipl}
\def\Oh{On the other hand}
 
\def\po{polynomial}

\def\sa{subalgebra}
\def\stf{straightforward}

\font\double=msbm10
\def\BbbC{\hbox{\double\char'103}}

\def\Dscr{\hbox{\double\char'104}}
\def\Tscr{\hbox{\double\char'124}}
\vbox{\vskip 4truecm}

\ctrline{\bf New rational solutions}
\ctrline{\bf of Yang-Baxter equation and deformed Yangians}
\ska
\centerline{\bf Alexander Stolin }
\centerline{\bf Department of Mathematics, Royal Institute of Technology }
\centerline{\bf S-100 44 Stockholm, Sweden}
\ske
\centerline{\bf Petr P. Kulish}
\centerline{\bf St.Petersburg Branch of Steklov Mathematical Institute}
\centerline{\bf  Fontanka 27, St.Petersburg 191011, Russia}
\ska
{\bf Abstract.}\ \ In this paper a class of
new quantum groups is presented: deformed Yangians. They arise from
rational solutions of the classical Yang-Baxter equation 
of the form ${c_2\over u}+const$.
The universal quantum R-matix for a deformed Yangian is
described. Its image in finite-dimensional representations
of the Yangian gives new matrix rational solutions of the 
Yang-Baxter equation (YBE).

\tf
{\bf 1. \qquad Introduction}
\par
The term ``quantum groups'' and the algebraic constructions 
associated with them appeared approximately 10 years ago in 
[D], [D2], [J1]. One of the starting points for such 
constructions was the classification of trigonometric
solutions of the classical Yang-Baxter equation (CYBE) 
obtained in [BD]. In particular $U_q (sl(n))$ can be viewed
as a quantization of the Lie bialgebra arising from the
Drinfeld-Jimbo solution of CYBE. Another ``quantum group''
was called Yangian ([D2]) and it arose from a rational 
solution of CYBE, a so-called Yang solution. Here we present
an attempt to define new quantum groups, which arise from
other rational solutions of CYBE.
\par
Now let $g$ be a simple \li\ over $\BbbC$.
\par
{\bf Definition.} \ \ Let $X(u,v)={c_2\over u-v}+r(u,v)$ be a function
from $\BbbC^2$ to $g\otimes g$. We say that $X(u,v)$ is a rational
solution of the classical Yang-Baxter equation (CYBE) if:
\par
\item{(1)}$c_2=\sum_\mu I_\mu\otimes I_\mu$, where $\{I_\mu\}$ is an
orthogonal basis of $g$ with respect to the Killing form;
\item{(2)}$r(u,v)$ is a \po\ in $u,v$ ;
\item{(3)}$X(u,v)=-X(v,u)^\s$, where $\s$ interchanges factors in $g
\otimes g$ ;
\item{(4)}$[X^{12}(u_1,u_2),\ X^{13}(u_1,u_3)] + 
[X^{12}(u_1,u_2),\ X^{23}(u_2,u_3)] + 
[X^{13}(u_1,u_3), \ X^{23}(u_2,u_3)]=0$. 
\par
Here $[X^{12},X^{13}]$ is the usual commutator in the associative
algebra $U(g)^{\otimes 3}$. The other two summands are defined in the
same way.
\ske
{\bf Definition.} \ \ We say that two rational solutions $X_1(u,v)$
and $X_2(u,v)$ are gauge equivalent if there exists an automorphism
$\lambda$ of algebra $g[u]$ such that $(\lambda\otimes\lambda)X_1(u,v)=
X_2(u,v)$.
\par
It turns out that the degree of the \po\ part of a rational solution
of CYBE can be estimated. More exactly, the following result was proved
in [S]:
\par
{\bf Theorem 1.} \ \ Let $X(u,v)=\ {c_2\over u-v}\ + r(u,v)$ be a rational
solution. Then there exists a rational solution $X_1(u,v)$, which is
gauge equivalent to $X(u,v)$ and such that
$$ X_1(u,v)=\ {c_2\over u-v}\ + a_0+b_1 u+b_2 v+cuv \ . $$
Here $a_0,b_1,b_2,c\in g^{\otimes 2}$.
\par
In the present paper we will be dealing with the case $X(u,v)=\ {c_2
\over u-v}\ + r_0$, where $r_0\in g^{\otimes 2}$. Clearly $X(u,v)$ is 
a solution of CYBE if and only if $r_0$ itself is a solution of CYBE.
\par
Let $K=\BbbC ((u^{-1}))$. One can define the following 
non-degenerate ad-invariant inner product on 
$g \otimes K$: \ \ $(x,y)=Res_{u=0} tr(adx \cdot ady).$ Denote
$g\otimes C[[u^{-1}]]$ by $g[[u^{-1}]]$.
\par
{\bf Theorem 2} (see [S]).\par
1) \ \ There is a 1-1 correspondence between the set of rational solutions
of CYBE of the form ${c_2\over u-v}\ + r_0$ and subalgebras $W \st g \otimes K$
such that: \par
({\bf i}) $u^{-2}g[[u^{-1}]] \st W \st  g[[u^{-1}]]$. \par
({\bf ii}) $W^{\bot}=W$ with respect to the form $(\ ,\ )$ introduced above.\par
({\bf iii}) $W \oplus g[u]=g\otimes K$.\par
2) \ \ Any $W$ satisfying conditions ({\bf i-iii}) above
defines a \sa\ $L\st g$ and a non-degenerate 2-cocycle $B$ on $L$. In
other words $B$ is skew-symmetric and satisfies
$$ B([x,y],z) + B([z,x],y)+B([y,z],x)=0 $$
for any $x,y,z\in L$. Moreover, $r_0$ is contained in $\Lambda^2 L$, is a non-degenerate
2-tensor and
$r_0^{-1}=B\in \Lambda^2 L^*$.
\par
A  Lie algebra with a non-degenerate 2-cocycle  is called quasiFrobenius. \par
3) \ \ Conversely, any pair $(L,B)$ such that $L$ is a \sa\ of $g$ and
$B$ is a non-degenerate 2-cocycle on $L$, defines a rational solution
of the form ${c_2\over u-v}\ + r_0$.
\par

Our approach to quantization of a rational solution of CYBE 
of the form $ {c_2\over u-v}\ + r_0$
is based on
the following result borrowed from [D1].
\par
{\bf Theorem 3.} \ \ Let $r_0 \in L\otimes L \st g\otimes g$ satisfy
CYBE. Then 
there exists an element $F\in (U(L)[[h]])^{\otimes 2}
\st (U(g)[[h]])^{\otimes 2}$ such that:
\par
1) \ \ $(\D_0\otimes 1)F\circ F^{12}=(1\otimes\D_0)F\circ F^{23}$, where
$\D_0:U(g)\ra U(g)^{\otimes 2}$ is the usual cocommutative co\mt ication.
\par
2) \ \ $F=1\otimes 1+\ \half\ hr_0+\sum^\nty_2 F_i h^i$ .
\par
3) \ \ $R=(F^{21})^{-1}F\in (U(L)[[h]])^{\otimes 2}
\st (U(g)[[h]])^{\otimes 2}$ satisfies YBE and is of the form
$R=1\otimes 1+\ \ hr_0+\sum^\nty_2 R_i h^i$ . Here $F^{21}=F^{\s},$
where $ {\s}$ interchanges factors in $ (U(g)[[h]])^{\otimes 2}$.


%
\ske
{\bf 2. \qquad Deformation of Yangians}
\par
Now we return to rational solutions of CYBE. The simpliest rational
solution is $X_0(u,v)=\ {c_2\over u-v}$, \ \ i.e.,\ \ $r(u,v)\equiv 0$.
Yangians were introduced by Drinfeld in [D2] in order to obtain a
``sophisticated quantization" of $X_0(u,v)$.
\par
{\bf Definition.} \ \ Let $g$ be a simple \li\ over $\BbbC$, given by
generators $\{I_\a\}$ and relations $[I_\a,I_\b]=C_{\a\b}^\g I_\g$,
where  $\{I_\g \}$ is an orthonormal basis with respect to the Killing form.
Then Yangian $Y(g)$ is an associative algebra with 1, generated by
elements $\{I_\a\}$ and $\{\Tscr_\a\}$ and the following relations
$$ [I_\a,I_\b]=C_{\a\b}^\g I_\g\ ;\ \ [I_\a,\Tscr_\b]=C^\g_{\a\b}
\Tscr_\g \ ; \ \leqno(1) $$
$$ [[\Tscr_\lambda[\Tscr_\mu,I_\nu]]-[I_\lambda[\Tscr_\mu,\Tscr_\nu]]=a_{
\lambda\mu\nu}^{\a\b\g}\{I_\a, I_\b, I_\g\} \ , \ \leqno(2) $$
here $a^{\a\b\g}_{\lambda\mu\nu}=\ {1\over 24} \ C^i_{\lambda\a}C^j_{\mu\b}C^k_{\nu
\g}C^k_{ij}$ and $\{x_1,x_2,x_3\}=\sum_{i\neq j\neq k}x_i x_j x_k$.
$$ \D I_\lambda=I_\lambda\otimes 1+1\otimes I_\lambda \ \leqno(3) $$
$$ \D\Tscr_\lambda=\Tscr_\lambda\otimes 1+1\otimes \Tscr_\lambda+\half\ C^\nu_{
\lambda\mu}I_\nu\otimes I_\mu \ . \  \leqno(4) $$
For any $a\in \BbbC$ define an automorphism $T_a$ of $Y(g)$ by formulas:
$$ T_a(I_\lambda)=I_\lambda \ ; \ \ T_a(\Tscr_\lambda)=\Tscr_\lambda+aI_\lambda \ . $$
As usual, we denote by $\D\pr$ the opposite co\mt ication.
\par
{\bf Theorem 4} ([D2]).
There exists a unique $R(u)= 1\otimes 1+\sum^\nty_{k=1}R_k u^{-k}, \
R_k\in Y(g)^{\otimes 2}$ such that

\item{1)}$(\D\otimes 1)R(u)=R^{13}(u)R^{23}(u)$ ;
\item{2)}$(T_u\otimes 1)\D\pr(a) = R(u)((T_u\otimes 1)\D(a))R(u)^{-1}$
\qquad for all $a\in Y(g)$ ;
\item{3)}$(T_a\otimes T_b)R(u)=R(u+a-b)$ ;
\item{4)}$R^{12}(u)R^{21}(-u)=1\otimes 1$ ;
\item{5)}$R^{12}(u_1-u_2)R^{13}(u_1-u_3)R^{23}(u_2-u_3)=
  R^{23}(u_2-u_3)R^{13}(u_1-u_3)R^{12}(u_1-u_2)$ ;
\item{6)}$R_1=c_2$ .
\hfill\black\par
The identity 2) means that $Y(g)$ is a pseudotriangular \Hf.
Consider $Y(g)[[h]]=Y_h(g)$. Clearly $Y_h(g)$
contains $U(g)[[h]]$ as a Hopf \sa. Let $F$ satisfy
condition 1 of Theorem 3 and we can view $F$ as an element 
of $(Y_h(g))^{\otimes 2}$.
Obviously, one can extend the \Hf\ structure to $Y_h(g)$.
Let us define a new algebra $\tilde Y_h(g)$, which has the same \mt
ication as $Y_h(g)$ but co\mt ication is defined as
$ \tilde\D(a)=F^{-1}\D(a)F \ . $
The main result of this paper is the following:
\par
{\bf Theorem 5.}
\item{1)}The algebra $\tilde Y_h(g)$ is a \Hf.
\item{2)}Define $\tilde R(u)$ to be $\tilde R(u)=(F^{21})^{-1}R(u)F$. \bk
Then $(\tilde\D\otimes 1)\tilde R(u)=\tilde R^{13}(u)\tilde R^{23}(u)$ ;
\item{3)}$(T_a\otimes T_b)\tilde R(u)=\tilde R(u+a-b)$ ;
\item{4)}$\tilde R^{12}(u)\tilde R^{21}(-u)=1\otimes 1$ ;
\item{5)}$(T_u\otimes 1)\tilde\D\pr(a)=\tilde R(u)((T_u\otimes 1)
\tilde\D(a))\tilde R(u)^{-1}$  for all $a\in \tilde Y_h(g)$ ;
\item{6)}$\tilde R^{12}(u_1-u_2)\tilde R^{13}(u_1-u_3)\tilde R^{23}(u_2-
u_3) = \tilde R^{23}(u_2-u_3)\tilde R^{13}(u_1-u_3)\tilde R^{12}
(u_1-u_2)$ ;
\item{7)}$\tilde R\left({u\over h}\right)=1\otimes 1 + h\left({
c_2\over u} +r\right) + 0(h)$ ;
\par
{\bf Proof.} \ \ 1) \ \ We must prove that $\tilde\D$ is a coassociative
operation. This is \stf\ from coassociativity $\D$ and the defining
identity for $F$.
\par
5) \ \ By the  
definition of $\tilde\D$ we have: $(T_u\otimes 1)\tilde\D\pr(a)=
(T_u\otimes 1)((F^{21})^{-1}\D\pr
(a)F^{21})$ .
\par
We note that $(T_a\otimes T_b)F=F$ since $F\in (U(g)[[
h]])^{\otimes 2}$. Hence,
$$ \eqalign{& (T_u\otimes 1)\tilde\D\pr(a) =  \cr
& = (F^{21})^{-1}R(u)((T_u\otimes 1)\D(a))R(u)^{-1}F^{21} \cr
& = ((F^{21})^{-1}R(u)F)((T_u\otimes 1)\tilde \D(a))((F^{21})^{-1}
R(u)F)^{-1} \cr
& = \tilde R(u)((T_u\otimes 1)\tilde\D(a)\tilde R(u)^{-1} \ \ {\rm
by\ Theorem\ 4.} \cr} $$
\par
2)  \ \ If $Y(g)$ were a triangular \Hf , all would follow
from results [D3]. It turns out that the pseudotriangular structure does
not affect considerations similar to ones of [D3]. We have:
$$ \eqalign{&(\tilde\D\otimes 1)\tilde R(u) =  \cr
& = (F^{12})^{-1}((\D\otimes 1)((F^{21})^{-1}R(u)F))F^{12} = \cr
& = (F^{12})^{-1}(\D\otimes 1)(F^{21})^{-1}((\D\otimes 1)R(u))(\D\otimes
1)F\circ F^{12} = \cr
& = (F^{12})^{-1}((\D\otimes 1)(F^{21})^{-1})\ (R^{13}(u)R^{23}(u))(1
\otimes\D)F\circ F^{23} \ .  \cr} $$
Again since $(T_a\otimes T_b)F=F$, it follows from Theorem 4 that
$ R^{23}(u)((1\otimes \D)F)=((1\otimes\D\pr)F)R^{23}(u) \ . $
\Oh\ $(1\otimes\D\pr)F=((\D\otimes 1)F)^{32}F^{13}(F^{32})^{-1}$, where
$(a\otimes b\otimes c)^{32}=a\otimes c\otimes b$.
Further, $R^{13}(u)((\D\otimes 1)F)^{32}=((\D\pr\otimes 1)F)^{32}R^{13}
(u)$. It remains to show, that
$$ (F^{12})^{-1}((\D\otimes 1)(F^{21})^{-1})\ ((\D\pr\otimes 1)F)^{32}=
(F^{31})^{-1} $$
which is true
by the defining relation for $F$.\par
3), 4) and 7) are \stf\ from the \corr ing statements of Theorem 4.
Let us deduce 6). It follows from 2) that
$ (T_a\otimes 1\otimes 1)((\tilde\D\otimes 1)\tilde R(x))=\tilde R^
{13}(x+a)\tilde R^{23}(x) \ . $
Hence,
$$ \tilde R(a)(T_a\otimes 1\otimes 1)\ ((\tilde\D\otimes 1)\tilde R(x))=
\tilde R^{12}(a)\tilde R^{13}(x+a)\tilde R^{23}(x) \ . $$
\Oh\ 5) implies that
$$ \tilde R(a)(T_a\otimes 1\otimes 1)\ ((\tilde\D\otimes 1)\tilde R(x))=
(T_a\otimes 1\otimes 1)\ ((\tilde\D\pr\otimes 1)\tilde R(x))\tilde R(a)
\ . $$
Since $(\tilde\D\pr\otimes 1)\tilde R(x)=\tilde R^{23}(x)\tilde
R^{13}(x)$,
we find that $(T_a\otimes 1\otimes 1)\ ((\tilde \D\otimes 1)\tilde R(x))
= \tilde R^{23}(x)\tilde R^{13}(x+a)$, which completes the proof.
\hfill\black\par

The results from [KST] show that the problem of finding 
explicit formulas leads to rather difficult computations
even in the simplest case of $s\ell(2)$ with $F$ found in [CGG].
Our aim is to present a number of cases when
a rational R-matrix for $s\ell(n)$ and $o(n)$ can be computed explicitly 
in the corresponding fundamental n-dimensional representations.
We need the following corrolary to Theorem 5:
\par
{\bf Corrolary 1.} Let ${c_2\over u-v}+r_0 $ be a rational solution
of CYBE for   $s\ell(n)$, $F\in U(s\ell(n))^{\otimes 2}$ be the
corresponding ``quantizing element'' and $R\in Mat(n,\BbbC)^{\otimes 2}$
be the image of the quantum R-matrix $(F^{21})^{-1}F$ in the 
fundamental n-dimensional representation of  $s\ell(n)$. 
If $P\in Mat(n,\BbbC)^{\otimes 2}$ is the permutation
matrix, which acts in ${\BbbC}^{n}\otimes {\BbbC}^n$ as $P(a\otimes b)=
b\otimes a$, then $uR+{P}$ satisfies YBE.\par
{\bf Proof.}\ \ Let us consider the R-matrix $\tilde R(u)=(F^{21})^{-1}R(u)F$, 
where $R(u)$ is Drinfeld's R-matrix for $Y(s\ell(n))$. It was proved in [D2] 
that
the image of  $R(u)$ in the n-dimensional representation is 
$1\otimes 1 +{P\over u}$ up
to a scalar factor. It 
is easy to see that $(T^{21})^{-1}PT=P$ for any invertible
$T\in Mat(n,\BbbC)^{\otimes 2}$. This observation completes the proof.\par

{\bf Remark.}\ \ It is worth noticing that we have proved that if 
$R\in Mat(n,\BbbC)^{\otimes 2}$
satisfies YBE and 
is unitary, i.e. 
$R^{21}R=1\otimes 1$, then $R+{P\over u}$ is a rational 
solution of YBE because
according to [D1]
any such $R$ comes from some $F\in U(g\ell(n))^{\otimes 2}$.
Of course, knowing the answer it is not difficult to check that
$R+{P\over u}$ is 
really a solution of YBE (using the fact that $(a\otimes b)P=
P(b\otimes a)$ for any $a,b\in Mat(n, \BbbC)$). However the general
approach provides rational solutions in any finite-dimensional
representation of $Y(s\ell(n))$. 
\par
\ske
{\bf Example 1.}\ \ We would like to expose a number of unitary 
R-matrices not involving complicated computations. According to
Theorem 2 we have to indicate a pair $(L,B)$, where $L\st s\ell(n)$
and $B$ is the corresponding non-degenerate 2-cocycle. Put
$$L=\{(a_{ij}): a_{ij}=0\ for\ i>j;\ a_{ii}=-a_{n+1-i,n+1-i}\},\ 
B(x,y)=f([x,y])$$
where $f(\{a_{ij}\})=\sum_{i+j=n+1}a_{ij}$.\par
Let $E_{ij}\in Mat(n)$ be the set of matrix units.
Let us denote $E_{ii}-E_{n+1-i,n+1-i}$ by $H_i$. Then the corresponding
classical r-matrix 
$r\in Mat(n)^{\otimes 2}$ has of the following form:

$$r={1\over 2}(\sum_i (H_i\otimes E_{i,n+1-i}- E_{i,n+1-i}\otimes
H_i))+\sum_{i<j<n+1-i} (E_{ij}\otimes E_{j,n+1-i}-E_{j,n+1-i}\otimes E_{ij})$$
Direct computations show that $r^3=0$.
It is known (see [CGG]) that in this case $R=1\otimes 1+r+{1\over 2}r^2$
is a unitary solution of YBE. Corrolary 1 implies that 
$$u(1\otimes 1 +r+{1\over 2}r^2) + P$$
is a rational solution of YBE.\par
\hfill\black
\ske
Let $o(N,\BbbC)$ be an orthogonal Lie algebra consisting of all matrices
$A\in Mat(N,\BbbC)$ such that $A^t=-A$.
Let $K\in  Mat(N,\BbbC)^{\otimes 2}$ be the matrix obtained from 
$P\in  Mat(N,\BbbC)^{\otimes 2}$ by the transposition in the first
factor. It was proved in [D2, KS] that the image of $R(u)\in Y(o(N))^
{\otimes 2}$ in the N-dimensional representation of $Y(o(N))$ is,
up to a scalar factor 
$$u1\otimes 1+ P-{u\over k+u}K,\ \ k={1\over 2}(N-2)$$
{\bf Corrolary 2.}\ \ Let $r\in o(N)^{\otimes 2}$ be a classical
r-matrix and $F\in U( o(N))^{\otimes 2}$ be the corresponding
quantizing element. Let us denote by $F_0$ 
(respectively $R_0$) the image of 
$F$ (respectively 
$(F^{21})^{-1}F$) in $Mat(N)^{\otimes 2}$.
\par
Then $uR_0+P-{u\over k+u}(F_0^{21})^{-1}KF_0$ satisfies YBE.\par
{\bf Proof.}\ \ The statement can be proved exactly as Corollary 1.
\par
{\bf Example 2.}\ \ Now we need an another realization of $o(N,\BbbC)$,
namely $$o(N)=\{(a_{ij})\in Mat(N):\ a_{ij}=-a_{N+1-j,N+1-i}\}$$
Let $T$ be any element of $GL(N,\BbbC)$ which conjugates the first
form of $o(N)$ to the second one. Denote $E_{11}-E_{NN}$ by $H$
and  $E_{12}-E_{N-1,N}$ by $E$. Clearly $e=TET^{-1}$ and $h=THT^{-1}$
are skew-symmetric matrices. Further we have $[h,e]=e$ since $[H,E]=E$ and
$r_0=h\otimes e - e\otimes h$ satisfies CYBE (for $N>3$).
The corresponding quantizing element $F$ was found in [CGG] and 
is of the form:
$$F=1\otimes 1 +\sum_n {1\over n!}h(h+1)...(h+n-1)\otimes e^n
\st U(o(N))^{\otimes 2}$$
Clearly $r_1=H\otimes E- E\otimes H,\ \ N>3$
also satisfies CYBE and therefore, we can compute the corresponding
matrix solution of YBE, which is 
$$1\otimes 1+r_1-E_{N-1,N}\otimes E_{12} \st 
Mat(N)^{\otimes 2}$$
since
the image of $E^2$ is $0$ in $Mat(N)$. Finally we obtain that the
following element of $Mat(N)^{\otimes 2}$ is a rational solution
of YBE:
$$(1\otimes 1+r_0-e_{-}\otimes e_{+})u+P-{u\over k+u}(1\otimes 1-
e\otimes h)K(1\otimes 1+
h\otimes e),\ \ k={1\over 2}(N-2)
$$
where $e_{-}=TE_{N-1,N}T^{-1}$ and  $e_{+}=TE_{12}T^{-1}$.\par
\ske
{\bf Remark.}\ \ It is interesting to point out that the construction of
the new 
solutions to the YBE (Theorem 5) preserves the regularity 
property [KS] of the initial $R$-matrix: $ R(0) = P $. Therefore 
one can obtain series of integrable models with local Hamiltonians 
corresponding to these new $R$-matrices. In the simplest case of $Y(s\ell(2))$ 
with non-standard quantization of $s\ell(2)$ (see [KST])
the spin-1/2 analog of the 
$XXX$-model on one-dimensional chain is given by the Hamiltonian 
$$     
H=\sum_n ( (\s_n, \s_{n+1}) + 
\xi^2 \s_n^- \s_{n+1}^- +  \xi (\s_n^- - \s_{n+1}^-) ), 
$$ 
where $\xi$ is a deformation parameter,  
$\s_n^x, \s_n^y, \s_n^z$ are Pauli sigma-matrices 
acting in $\BbbC_n^2$ related to the $n$-th site of the chain
and $\s_n^{-}={1\over 2}(\s_n^x-i\s_n^y)$.
\ske
{\bf Acknowledgements.}
The authors are thankful to Professors J. Lukiersky and A. Molev 
for valuable discussions. The results of this paper were delivered  
by the first author at the Colloquium
``Quantum Groups and Integrable Systems'' in Prague, June 1996.
The visit was supported by Swedish Natural Science Research Council.

\ska
\ctrline{\bf References}
\par
\item{[BD]\ \ \ \ \ }Belavin, A.A., Drinfeld, V.G., Funct. Anal. Appl., 
{\bf 16} (1982) 159.
\item{[CGG]\ \ }Coll, V., Gerstenhaber, M., Giaquinto, A., Israel Math. 
Conf. Proc., vol. {\bf 1}, Weizmann Science Press, (1989).
\item{[D]\ \ \ \ \ \ }Drinfeld, V.G., Proc. ICM-86 (Berkeley), 
vol. {\bf 1} (1986) 798.
\item{[D1]\ \ \ \ \ }Drinfeld, V.G., Soviet Math. Dokl., {\bf 27} (1983) 68.
\item{[D2]\ \ \ \ \ }Drinfeld, V.G., Soviet Math. Dokl., {\bf 32} (1985) 254.
\item{[D3]\ \ \ \ \ }Drinfeld, V.G., Leningrad Math. J., {\bf 1} (1990) 1459.
\item{[J1]\ \ \ \ \ }Jimbo, M., Lett. Math. Phys., {\bf 10} (1985) 63.
\item{[KS]\ \ \ \ \ }Kulish, P.P., Sklyanin, E.K., Zap. Nauch. Semin. 
LOMI, {\bf 95} (1980) 129.
\item{[KST]\ \ \ \ \ }Khoroshkin, S., Stolin, A., Tolstoy, V., 
preprint TRITA-MAT-1995-MA-17, KTH, Stockholm (1995).  
\item{[S]\ \ \ \ \ \ }Stolin, A., Math. Scand., {\bf 69} (1991) 56.
\bye